\documentclass[10pt]{amsart}

\usepackage[T1]{fontenc}
\usepackage{amsmath,amssymb,mathtools}
\usepackage{mathrsfs}
\usepackage{microtype}
\usepackage{tikz}
\usepackage{tikz-cd}
\usepackage{subcaption}
\usepackage{booktabs}
\usepackage{tabularx}
\usepackage{longtable}
\usepackage{enumitem}
\usepackage{listings}
\usepackage{upquote}
\usepackage[margin=0.95in]{geometry}
\usepackage[colorlinks=true,citecolor=blue,linkcolor=blue,urlcolor=blue]{hyperref}

\newtheorem{theorem}{Theorem}[section]
\newtheorem{proposition}[theorem]{Proposition}
\newtheorem{lemma}[theorem]{Lemma}

\theoremstyle{definition}

\theoremstyle{remark}

\theoremstyle{plain}
\newtheorem*{theoremA}{Theorem A}
\newtheorem*{theoremB}{Theorem B}

\title[Optimality of the $D_4$ lattice quantizer]{Optimality of the $D_4$ Lattice Quantizer}
\author{Tianle Liu}
\address{University of Southern California}
\email{tianleli@usc.edu}
\date{}

\makeatletter
\@ifundefined{subjclassname@2020}{%
  \@namedef{subjclassname@2020}{%
    \textup{2020} Mathematics Subject Classification}%
}{}
\makeatother
\hypersetup{
  pdfauthor={Tianle Liu},
  pdftitle={Optimality of the D4 Lattice Quantizer}
}

\begin{document}

\begin{abstract}
We prove that the root lattice $D_4$ uniquely minimizes the normalized
second moment among four-dimensional lattices, up to similarity.
The minimum is $13/(120\sqrt{2})$.
The proof reduces the problem to polynomial inequalities on three
secondary cones. We derive the second-moment formulas
geometrically and give explicit polynomial certificates that can be
checked by coefficient comparison.
We also derive a lower bound for lattice quantizers in every dimension
at least two from the number of Voronoi facets, using a comparison with
circular cones and concavity in the solid angle.
This bound strictly improves the sphere bound and recovers the
hexagonal value in dimension two.
\end{abstract}
\maketitle

\section{Introduction}\label{sec:introduction}

Lattice quantization represents a real-valued vector by its nearest
lattice point. At a fixed point density, the normalized second moment
of the Voronoi cell determines the mean squared quantization error
under a uniform cell distribution. The lattice quantizer problem asks
which lattice minimizes this quantity in each dimension.

For a full-rank lattice $L\subset\mathbb R^d$, let
\[
 V(L)=\{x\in\mathbb R^d:\ \|x\|\le\|x-\lambda\|
                         \text{ for every }\lambda\in L\}
\]
be its Voronoi cell at the origin, and put $v=\operatorname{covol}(L)$.
The normalized second moment is
\begin{equation}\label{eq:def-G}
 G(L)=\frac{1}{d v^{1+2/d}}\int_{V(L)}\|x\|^2\,dx.
\end{equation}
It is invariant under similarities. In particular, minimizing $G$
is equivalent to minimizing the second moment at a fixed cell volume.

\subsection{The optimal lattice in dimension four}

Write
\[
 D_4=\{z\in\mathbb Z^4:z_1+z_2+z_3+z_4\in2\mathbb Z\}.
\]
\begin{theoremA}
Every full-rank lattice $L\subset\mathbb R^4$ satisfies
\[
 G(L)\ge G(D_4)=\frac{13}{120\sqrt2}.
\]
Equality holds if and only if $L$ is similar to $D_4$.
\end{theoremA}

For a uniform planar source, Fejes T\'oth established the asymptotic
optimality of the hexagonal arrangement \cite{FejesToth1959}.
Gersho proposed the body-centered cubic structure in dimension three
and the structure associated with $D_4$ in dimension four
\cite{Gersho1979}. Conway and Sloane computed the Voronoi geometry and
exact normalized second moments of classical root lattices, including
the regular $24$-cell of $D_4$ \cite{CS82}, and gave efficient
nearest-point algorithms \cite{CSFast1982}. Barnes and Sloane proved
that the body-centered cubic lattice is the unique optimum among
three-dimensional lattices, up to similarity \cite{BS83}.
Vel\'asquez Soto derived explicit moment formulas for the three
four-dimensional Voronoi types and proved local minimality of $D_4^*$,
which is similar to $D_4$
\cite[\S4.5 and Theorem~4.7.2]{VelasquezSoto2012}.
Zimmermann proved local minimality of $D_4$
\cite[Theorem~6.4.3]{Zimmermann2017}; recent numerical optimization and
exact identification provide further evidence for its optimality
\cite{APA25}. Theorem~A proves that $D_4$ is the unique global minimizer
among four-dimensional lattices, up to similarity.

For Theorem~A, Voronoi's reduction theory gives three full-dimensional
secondary cones, up to integral change of basis
\cite[\S6.4.3]{Zimmermann2017}. Their closures cover the positive-definite
forms. On each cone, integration expresses the second moment in
terms of the determinant and a homogeneous polynomial of degree five.
Six polynomial identities express suitable differences using squares,
Schur cubics, and polynomials with nonnegative coefficients. They reduce
the comparison to one-variable inequalities, give a strict inequality
on the principal cone, and identify a single equality ray on each
remaining cone.
Both rays represent $D_4$, and continuity includes all positive-definite
boundary forms. Section~\ref{sec:d4-global} gives the geometric moment
formulas and equality argument. All polynomial certificates are given
in Appendix~\ref{app:explicit-certificates}; their verification consists
of finite coefficient comparisons of the displayed expressions.
Computational searches were used to find these certificates, but no
external coefficient files are needed to check them.

\subsection{Bounds from Voronoi facets}

In dimensions where the optimum is unknown, the classical sphere bound
\begin{equation}\label{eq:sphere-bound}
 G(L)\ge\frac{1}{(d+2)\kappa_d^{2/d}},
 \qquad \kappa_d=\frac{\pi^{d/2}}{\Gamma(1+d/2)},
\end{equation}
follows by replacing the Voronoi cell with a ball of the same volume
\cite{Zador1982,CS82}. A second result improves this bound using the
number of facets. Put $A_d=d\kappa_d$, the surface area of the unit
sphere, and let $g_d(\omega)$ be the normalized second moment, about
the vertex, of a right circular cone of solid angle
$0<\omega<A_d/2$. Its explicit formula is given in
Section~\ref{sec:cone-bound}.

\begin{theoremB}
Let $d\ge2$, and let $L\subset\mathbb R^d$ be a full-rank lattice
whose Voronoi cell has $m$ facets. Put $f_d=2(2^d-1)$. Then
\[
 G(L)\ge m^{-2/d}g_d(A_d/m)
       \ge f_d^{-2/d}g_d(A_d/f_d).
\]
\end{theoremB}

The first inequality uses the actual facet count; the second uses the
universal estimate $m\le f_d$. Both bounds are strictly stronger than
\eqref{eq:sphere-bound}, and the dimension-dependent bound recovers
the hexagonal value in dimension two. In dimension four, the cone
comparison is strict even for $D_4$ and does not give Theorem~A.
The proof of Theorem~B is independent of the four-dimensional
classification: it compares each facet pyramid with a circular cone
and combines the comparisons by convexity and concavity.
It proves the lattice case of Ajjanagadde's proposed quantization
bound and his circular-cone concavity conjecture
\cite[Chapter~4, Conjectures~2 and~4]{Ajjanagadde2020}.

\subsection{Distortion bounds for lattice quantizers}\label{sec:distortion}

Let $L\subset\mathbb R^d$ be a full-rank lattice with cell volume $v$,
and let $Q_L$ denote Euclidean nearest-point quantization. Let $U$ be
uniform on $V(L)$, independent of the input $X$, and available to both
encoder and decoder. For the infinite lattice quantizer without
overload, subtractive dither gives
$\widehat X=Q_L(X+U)-U$. The error $\widehat X-X$ is uniform on $V(L)$
and independent of $X$ \cite{ZamirFeder1996}, so
\[
 \frac1d\mathbb E\|\widehat X-X\|^2=G(L)v^{2/d}.
\]

At fixed cell volume, Theorem~A therefore determines the minimum
mean squared error among four-dimensional lattice quantizers in this
model. Compared with independent uniform scalar quantization with equal
step sizes, whose lattice is a scaled copy of $\mathbb Z^4$, the optimal
distortion ratio is
\[
 \frac{G(D_4)}{G(\mathbb Z^4)}=\frac{13}{10\sqrt2}.
\]
Nearest-point quantization in $D_4$ requires coordinate rounding and,
if the sum is odd, changing the coordinate with largest absolute
rounding residual to the nearest alternative integer \cite{CSFast1982}.

If $V(L)$ has $m$ facets, Theorem~B gives
\[
 \frac1d\mathbb E\|\widehat X-X\|^2
 \ge v^{2/d}m^{-2/d}g_d(A_d/m).
\]
It also gives, for any candidate
lattice $L_0$ of dimension $d\ge2$,
\[
 f_d^{-2/d}g_d(A_d/f_d)
 \le\inf_{\dim L=d}G(L)\le G(L_0).
\]
The difference between these bounds limits the possible improvement
over $L_0$ at the same cell volume.

For orthogonal blocks of equal dimension, let
$L_1,\ldots,L_k\subset\mathbb R^d$ be full-rank lattices with cell
volumes $v_i$. Their Voronoi cells form a Cartesian product, and
\[
 G(L_1\times\cdots\times L_k)
 =\frac{\frac1k\sum_i G(L_i)v_i^{2/d}}
        {(\prod_i v_i)^{2/(kd)}}.
\]
A common lower bound $G(L_i)\ge b>0$ therefore implies
$G(L_1\times\cdots\times L_k)\ge b$ by the arithmetic--geometric
mean inequality. Equality holds exactly when every $G(L_i)=b$ and
$v_1=\cdots=v_k$. For four-dimensional blocks, Theorem~A gives the sharp
value $b=G(D_4)$, attained precisely by equally scaled $D_4$ blocks,
up to independent orthogonal transformations.

\subsection{Further background}

The signal-processing origins of quantization include Bennett's analysis
of quantization noise in pulse-code modulation \cite{Bennett1948} and
Lloyd's least-squares design conditions, developed in 1957 and published
in 1982 \cite{Lloyd1982}. In high-resolution theory, Zador separated the
source-dependent contribution to asymptotic distortion from a constant
depending on dimension and the distortion measure \cite{Zador1982}.
The defining optimization ranges over all quantizers.
Gray and Neuhoff \cite{GrayNeuhoff1998} give a detailed
account of these developments and of the distinction between
finite-rate design and high-resolution analysis.

Ajjanagadde's proposed facet bound extends the planar method by
replacing the triangles over polygon edges with pyramids over Voronoi
facets \cite[\S4.3.4]{Ajjanagadde2020}. Section~\ref{sec:cone-bound}
gives the required pyramid comparison and the precise relation to
his conjectures. His improvement of the sphere bound for unrestricted
quantizers in sufficiently large dimensions instead uses truncated
balls and variation in cell volumes
\cite[Theorem~10]{Ajjanagadde2020}.

\subsection{Quantization in learning, compression, and retrieval}

Computing systems use quantization to reduce storage, arithmetic cost,
and communication for model weights, activations, gradients, and
learned data representations. A direct application of lattice
quantization is UVeQFed, which uses subtractive dither for federated
model updates \cite{Shlezinger2021UVeQFed}. At a fixed cell volume,
the normalized second moment determines the mean squared error of
the dithered lattice quantizer, as described in Section~\ref{sec:distortion}.
For randomized gradient quantization, QSGD analyzes the tradeoff
between communication and the variance introduced by quantization
\cite{Alistarh2017QSGD}.

For model weights, AQLM represents a block by a sum of entries from
learned codebooks \cite{egiazarian2024aqlm}.
Lattice codebooks give such vector representations an explicit
algebraic structure. QuIP\# combines randomized Hadamard transforms
with a finite codebook derived from $E_8$ and also compares $D_4$
codebooks \cite[\S4 and Appendix~C.4]{tseng2024quipsharp}.
The generating matrix can also be learned separately for each weight
group, as in grouped lattice vector quantization \cite{Zhang2025Grouped}.
Kohli et al. study $D_4$ and $E_8$ in low-bit quantization-aware
training \cite{kohli2026lattice}, while E8-QAT uses an implicit
$E_8$ codebook during training \cite{Lei2026E8QAT}.
Higher-dimensional constructions include Leech-lattice codebooks
for language-model compression \cite{Ouderaa2026LLVQ}.
Related structured approaches include QTIP, which uses trellis-coded
quantization \cite{Tseng2024QTIP}, and pyramid vector quantization,
which represents weight directions by a fixed integer construction
\cite{Ouderaa2024PVQ}.

Learned media compression uses quantization in the latent space of
an encoder. EnCodec uses residual vector quantization for audio
\cite{Defossez2023EnCodec}. For images, Zhang and Wu learn lattice
generating matrices adapted to the latent distribution
\cite{Zhang2024OptimalLattice}; Lei, Hassani, and Saeedi Bidokhti
study lattice transform coding, including block and nested-lattice
constructions \cite{Lei2025LatticeTransform}.
Spherical Leech quantization uses a structured lattice-derived
codebook for visual tokenization and generation \cite{Zhao2025Leech}.

For neural-network inference, Deep Compression combines pruning,
learned weight sharing, and entropy coding \cite{han2016deep}, while
Jacob et al. develop training procedures for integer-only inference
with quantized weights and activations \cite{Jacob2018Integer}.
For large language models, LLM.int8() treats activation outliers in
higher precision \cite{Dettmers2022LLMint8}, while SmoothQuant reduces
them by rescaling weights and activations \cite{Xiao2023SmoothQuant}.
GPTQ, published as OPTQ, uses approximate second-order information to
quantize weights \cite{frantar2023optq}; AWQ uses activation statistics
to choose weight-scaling factors \cite{Lin2024AWQ}.

Quantization also reduces the memory occupied by attention keys and
values. KIVI uses different quantization groupings for the key and
value caches \cite{Liu2024KIVI}. TurboQuant combines random rotations
with scalar quantization and, for inner-product estimation, a correction
of the residual error; its applications include attention caches and
nearest-neighbor search \cite{Zandieh2026TurboQuant}.

Training raises additional questions about numerical range and the
propagation of rounding errors. Stochastic rounding has been used
in fixed-point training \cite{Gupta2015Limited} and in training
with eight-bit floating-point numbers \cite{Wang2018EightBit}.
The FP8 formats studied by Micikevicius et al. allocate exponent
and mantissa bits to accommodate different numerical ranges
\cite{Micikevicius2022FP8}. Wang et al. reduce the cost of generating
random numbers in integer training by sharing stochastic-rounding
thresholds and correcting the resulting structured rounding error
\cite{Wang2026StructuredRandomness}.
For model adaptation, QLoRA keeps a four-bit base model
fixed and trains low-rank adapters \cite{Dettmers2023QLoRA}.

In retrieval, product quantization partitions an embedding into
blocks and represents each block by a codeword \cite{Jegou2011}.
Anisotropic vector quantization adapts the distortion criterion to
inner-product search \cite{Guo2020}. Recommendation models use
quantized embeddings to allocate storage across categorical features
\cite{Kang2020}, and TIGER uses residual vector quantization to
construct item identifiers for generative retrieval \cite{Rajput2023}.

These applications involve source distributions, finite codebooks,
clipping, and task-dependent losses. Theorems~A and~B provide geometric
benchmarks for lattice blocks under the Euclidean distortion model
of Section~\ref{sec:distortion}. Stochastic rounding additionally
requires controlling randomization and error propagation.

\subsection{Notation}

Throughout the paper, the columns of a basis matrix $B$ generate
$B\mathbb Z^d$, and its Gram matrix is $Q=B^TB$. We also write
$G(Q)=G(B\mathbb Z^d)$. In these coordinates,
\begin{equation}\label{eq:gram-error}
 G(Q)=\frac{1}{d(\det Q)^{1/d}}\int_{[0,1]^d}\min_{z\in\mathbb Z^d}
                       (x-z)^TQ(x-z)\,dx.
\end{equation}
All integrals use Lebesgue measure, except solid angles, which use
surface measure on the unit sphere.

\section{The optimal four-dimensional lattice quantizer}\label{sec:d4-global}

\begin{theorem}\label{thm:d4-global}
Every full-rank lattice $\Lambda\subset\mathbb R^4$ satisfies
\[
 G(\Lambda)\ge\frac{13}{120\sqrt2}.
\]
Equality holds if and only if $\Lambda$ is similar to $D_4$.
\end{theorem}

The proof uses the three secondary cones in dimension four. On each
cone, integration expresses the normalized second moment in terms
of the determinant and a homogeneous polynomial of degree five.
Six polynomial identities give lower bounds for these polynomials and
reduce the comparison to one variable.

Voronoi's second reduction theory gives three full-dimensional secondary
cones in dimension four, up to integral change of basis. We use the
generators recorded in \cite[\S6.4.3]{Zimmermann2017}. Set
\[
 E_{0i}=e_ie_i^T,\qquad E_{ij}=(e_i-e_j)(e_i-e_j)^T\quad(1\le i<j\le4),
\]
and
\[
 A=\begin{pmatrix}
 2&1&-1&-1\\1&2&-1&-1\\-1&-1&2&0\\-1&-1&0&2
 \end{pmatrix},\qquad
 E_b=(1,1,-1,-1)^T(1,1,-1,-1).
\]
The closed cones $C_1,C_2,C_3$ are generated, respectively, by all ten
$E_{ij}$; by the nine $E_{ij}$ other than $E_{12}$ together with $A$; and
by the eight $E_{ij}$ other than $E_{12},E_{34}$ together with $E_b,A$.
After an integral change of basis, every positive-definite form belongs
to one of these cones.

For $Q$ positive definite put $D=\det Q$ and
\[
 M(Q)=\int_{[0,1]^4}\min_{z\in\mathbb Z^4}(x-z)^TQ(x-z)\,dx,
 \qquad G(Q)=\frac{M(Q)}{4D^{1/4}}.
\]
Define $N_i=60DM(Q)$ on $C_i$, for $i=2,3$, and $S=12DM(Q)$ on $C_1$.
We compute these homogeneous degree-five polynomials directly from the
Voronoi cells.

The moment $M$ is continuous on the positive-definite cone. Indeed, in
a compact neighborhood of a positive-definite form the eigenvalues are
bounded above and away from zero. For $x\in[0,1]^4$, comparison with
$z=0$ bounds the distance of every minimizing $z$ uniformly. Thus the
minimum is taken over a common finite set, and integration preserves
continuity. Formulas established for positive cone parameters therefore
extend to all boundary parameters for which $Q$ remains positive definite.

\subsection{The moment polynomial on the principal cone}

Index the edges of $K_5$ by $0\le i<j\le4$, and put
$g_{0i}=-e_i$ and $g_{ij}=e_i-e_j$ for $1\le i<j\le4$.
The sign convention on $g_{0i}$ is immaterial.  Write
\[
 Q=\sum_e\lambda_e g_eg_e^T,\qquad
 D=\det Q,\qquad L=\sum_e\lambda_e,
 \qquad \lambda_e>0.
\]
Every connected spanning graph $H$ with five edges on
$\{0,1,2,3,4\}$ has a unique cycle; write $\ell(H)$ for its length.

\begin{proposition}
The quantity $S=12D M(Q)$ is
\[
 S=DL+\sum_H\bigl(\ell(H)-2\bigr)\prod_{e\in H}\lambda_e,
\]
where $H$ ranges over these connected spanning graphs.
There are $150$, $60$, and $12$ such graphs with cycle lengths
$3$, $4$, and $5$, respectively. In particular,
$S\ge DL$ throughout the positive-definite part of the closed cone.
\end{proposition}

\begin{proof}
Use the dual coordinates $r=2Qx$.  The Voronoi inequalities become
$z^Tr\le z^TQz$ for $z\in\mathbb Z^4$.  Consider the zonotope
\[
 Z=\sum_e[-\lambda_e g_e,\lambda_e g_e].
\]
For every integer $z$ its support function satisfies
\[
 h_Z(z)=\sum_e\lambda_e|z^Tg_e|
       \le\sum_e\lambda_e(z^Tg_e)^2=z^TQz.
\]
Here the inequality follows from $|m|\le m^2$ for integers $m$.
Thus $Z$ is contained in the Voronoi polytope in $r$-coordinates.
The matrix $G=(g_e)_e$ is a reduced incidence matrix.  Its nonzero
$4\times4$ minors are $\pm1$ and correspond exactly to the $125$
spanning trees of $K_5$.  The zonotope volume formula and Cauchy--Binet
therefore give
\[
 \operatorname{vol}Z=16\sum_T\prod_{e\in T}\lambda_e=16D.
\]
The Voronoi polytope has the same volume, because the Jacobian of
$x\mapsto2Qx$ is $16D$.  The containment is consequently an equality.

We integrate $Z$ using a tiling by parallelotopes.
Order the ten edges lexicographically and
choose lifting heights $h_e=2^{\operatorname{index}(e)}$, with indices
starting at zero.  For a spanning tree $T$, write $G_T$ for its four
columns and set
\[
 \epsilon_{T,e}=\operatorname{sgn}
       (h_e-h_T^TG_T^{-1}g_e),\qquad
 c_T=\sum_{e\notin T}\epsilon_{T,e}\lambda_e g_e.
\]
The arguments of these signs are nonzero: the fundamental circuit of
$e$ has coefficients in $\{0,\pm1\}$, and its largest lifting height
exceeds the sum of all preceding heights. The upper facets of the
lifted zonotope $\sum_e[-\lambda_e(g_e,h_e),\lambda_e(g_e,h_e)]$ project to a tiling
of $Z$ by
\[
 P_T=c_T+\sum_{e\in T}[-\lambda_e g_e,\lambda_e g_e].
\]
Indeed, an upper supporting normal can be normalized to have final
coordinate one; the displayed signs specify the fixed endpoints of
all generators outside its four-dimensional face.  The upper boundary
is the graph of a function over $Z$, so these projections cover $Z$
and have disjoint interiors.

Put $K=\operatorname{adj}Q$ and $\lambda_T=\prod_{e\in T}\lambda_e$.
Since $\operatorname{vol}P_T=16\lambda_T$, integration over each
parallelotope gives
\[
 \int_{P_T}r^TKr\,dr
 =16\lambda_T\left(c_T^TKc_T+
              \frac13\sum_{e\in T}\lambda_e^2g_e^TKg_e\right).
\]
The change of variables above gives the polynomial identity
\begin{equation}\label{eq:c1-independent-integration}
 DS=\sum_T\lambda_T
       \left(3c_T^TKc_T+
                    \sum_{e\in T}\lambda_e^2g_e^TKg_e\right).
\end{equation}
Expanding the right side and collecting monomials in the $\lambda_e$
gives $D$ times the stated degree-five expression.
The identity extends by continuity to
nonnegative weights whenever $Q$ remains positive definite.
\end{proof}

This is the moment formula for the principal cone recorded in
\cite[\S6.4.3, equation~(6.10)]{Zimmermann2017}. The two stronger
inequalities used below are proved in
Proposition~\ref{prop:d4-exact-sos-certificates}.

\subsection{Moment formulas on the other two cones}

Let $\mathcal D(z)=\sum_Tz_T$ be the spanning-tree polynomial of $K_5$ and
\[
 \Phi(z)=\mathcal D(z)\sum_ez_e+\sum_H(\ell(H)-2)z_H,
 \qquad z_J=\prod_{e\in J}z_e,
\]
where $H$ ranges over connected spanning graphs with five edges and
$\ell(H)$ is the length of the unique cycle. Thus $S=12DM=\Phi$ on $C_1$.
For $C_2$ use the parameter order
\[
 (a,u,v,p,q,y_{13},y_{14},y_{23},y_{24},r),
\]
corresponding to
\[
 (a,\lambda_{01},\lambda_{02},\lambda_{03},\lambda_{04},
 \lambda_{13},\lambda_{14},\lambda_{23},\lambda_{24},\lambda_{34}).
\]
Set
\[
 z^{(2)}=(u+a,v+a,p,q,-a,y_{13}+a,y_{14}+a,y_{23}+a,y_{24}+a,r).
\]
On $C_3$ replace the final parameter $r$ by $b$, and put
\[
 \theta=(a+b,u-b,v-b,p,q,y_{13},y_{14},y_{23},y_{24},-b).
\]
\begin{proposition}\label{prop:d4-reconstructed-moments}
With $N_i=60DM$ on $C_i$, the moment and determinant formulas are
\begin{align*}
 D_2&=\mathcal D(z^{(2)}),&
 N_2&=5\Phi(z^{(2)})+
 2a^3\bigl[2a^2+5a(p+q+r)+10(pq+pr+qr)\bigr],\\
 D_3&=D_2(\theta),&
 N_3&=N_2(\theta)+2b^3\bigl[3b^2-5b(u+v)+10uv\bigr].
\end{align*}
Moreover, $N_i\ge20Dt$ on both cones, where
$t=a+\frac14\sum_j\lambda_j$, and $M(A)=13/30$.
\end{proposition}

The two correction terms come from changes in the diagonal of a
Delaunay square. We first identify the affected cells and reduce the
change in the moment to an area integral. The same integral formula
then gives both correction terms.

\paragraph{Delaunay triangulations}
The principal Delaunay triangulation consists, modulo $\mathbb Z^4$, of
the $24$ staircase simplices
\[
 [0,e_{\pi(1)},e_{\pi(1)}+e_{\pi(2)},
 e_{\pi(1)}+e_{\pi(2)}+e_{\pi(3)},e_1+e_2+e_3+e_4],\qquad \pi\in S_4.
\]
Write $c=e_i+e_j$. For each link edge $(V,W)$ below, replace
\[
 \mathcal T_-=
 \{[0,e_i,c,V,W],[0,e_j,c,V,W]\}
 \quad\hbox{by}\quad
 \mathcal T_+=
 \{[0,e_i,e_j,V,W],[e_i,e_j,c,V,W]\}.
\]
Both are triangulations of the same four-dimensional polytope, the
join of the square $[0,e_i,c,e_j]$ and the segment $[V,W]$.
The affine relation $0+c=e_i+e_j$ has quadratic lifting difference
$2q_{ij}$, so its sign selects the diagonal.

For the first change take $(i,j)=(1,2)$ and the six link edges
\[
\begin{array}{cc}
(-e_3,-e_3-e_4)&(-e_4,-e_3-e_4)\\
(-e_3,c+e_4)&(-e_4,c+e_3)\\
(c+e_3,c+e_3+e_4)&(c+e_4,c+e_3+e_4).
\end{array}
\]
These replace twelve of the staircase simplices. For the second change
take $(i,j)=(3,4)$ and the four link edges
\[
 (-e_2,c+e_1),\quad(-e_1,c+e_2),\quad
 (c+e_2,c+e_1),\quad(-e_1,-e_2).
\]
These replace eight of the simplices obtained after the first change.
In each case the old simplices occur without repetition modulo lattice
translation. Replacing triangulations of the same cells therefore
preserves the covering and disjoint interiors.
The $120$ facet occurrences of the $24$ simplices pair into $60$
translation classes, each containing exactly two facets with opposite
orientations. This follows by comparing the vertex differences in the
displayed simplex lists. The paired facets have identical induced
subdivisions, so the replacements agree on cell boundaries.

To check that the resulting triangulations are Delaunay, we compare
the affine interpolants across their common facets. Each has $24$ unimodular simplices and
$60$ facets modulo translation. At a common facet, let $L_\sigma$ be
the affine function interpolating $Q$ at the five vertices of one
simplex, and let $w$ be the opposite vertex of the other simplex.
The distinct forms $Q[w]-L_\sigma(w)$ and their multiplicities are
\[
\begin{array}{c|c|c}
\text{cone}&\text{lifting difference}&\text{multiplicity}\\\hline
C_1&2\lambda_e&6\text{ for each of ten edges }e\\
C_2,C_3&2a&6\\
C_2,C_3&2\lambda_j&4\text{ for each of nine parameters }\lambda_j\\
C_2,C_3&2(a+\lambda_j)&2\text{ for each of nine parameters }\lambda_j.
\end{array}
\]
This table is obtained directly from the displayed lattice vertices by
affine interpolation. Every entry is positive in the corresponding
open cone. The piecewise affine interpolant of $Q$ on the periodic
triangulation is consequently locally convex across every facet, hence
globally convex. Indeed, its restriction to a line avoiding the
codimension-two faces has nondecreasing slopes; passage to the limit
handles all lines. Its affine pieces are therefore supporting
hyperplanes of the lifted lattice, proving the Delaunay assertion.

\paragraph{Reduction to an area integral}
Write $Q[n]=n^TQn$ and use dual coordinates $r=2Qx$. Then
\[
 f_Q(r)=\frac14r^TQ^{-1}r+
          \min_{n\in\mathbb Z^4}\{Q[n]-n^Tr\},\qquad
 N(Q)=\frac{15}{4}\int_{\mathbb R^4/(2Q\mathbb Z^4)}f_Q(r)\,dr.
\]
The six, respectively four, cells listed above are the complete set of
cells that change at a positive point of the common boundary of the
corresponding cones. Their common square has a compact dual face $F$: a
hexagon between $C_1$ and $C_2$, and a rectangle between $C_2$ and $C_3$. The union of the nearest-site sets over
the closed face consists precisely of the four square vertices and the
listed link sites.
Positive definiteness bounds all potentially competing lattice sites
on a bounded neighborhood, so every omitted site is uniformly separated
from the minimum there.

Each face embeds in the period torus. Write $S$ for its square and
$V=\operatorname{span}(S-S)$. If $r,r'=r+2Qn$ both belong to $F$,
subtracting their equal-distance equations gives $n\perp_QV$.
The nearest-site set at $r'$ contains both $S$ and $S+n$. If these
intersect, then $n\in V$, so positive definiteness gives $n=0$.
Otherwise it contains at least eight sites, whereas the listed Delaunay cells give at most six nearest sites at
each point of $F$. Thus $F$ embeds. Compactness and positive
definiteness give a small embedded neighborhood, uniformly near the
boundary, so the integral over this neighborhood is counted once per period.

Choose a small embedded polyhedral neighborhood $W_Q$ of $F$.
Only the integral over $W_Q$ contributes to the difference of the two
analytic expressions for the moment on the adjacent cones. To account for the moving
periods, start with $P_Q=2Q(c+[0,1]^4)$ and remove its intersections
with all translates of $W_Q$.
Only finitely many translates meet $P_Q$, uniformly near the boundary.
Periodicity identifies the sum of their integrals with one integral
over $W_Q$. The remaining region has a finite subdivision by nearest-site
regions and boundaries of the translated neighborhoods. Choose the artificial cuts generically,
away from the faces where the subdivision changes. Its combinatorial type is fixed; its
vertices solve fixed selections of nonsingular linear systems with
analytic coefficients.
A fixed triangulation therefore makes the remaining quadratic integrals
analytic across the boundary. The quadratic part of $f_Q$ on $W_Q$ is also
analytic. Thus the desired correction is $15/4$ times the difference between
the two analytic continuations of the integral of the minimum of the
affine functions over $W_Q$.

The next lemma computes this difference. It applies to both changes of
triangulation after translating the two coordinates transverse to $F$.

\begin{lemma}\label{lem:d4-area-integral}
Let $K\subset\mathbb R^2$ be a fixed compact polygon. For $t$ near zero,
let $U_t,V_t\le0$ be continuous piecewise affine functions on $K$.
Suppose their common subdivision has fixed combinatorial type, its
vertices depend analytically on $t$, and the coefficients of the affine
functions on each cell depend analytically on $t$. Put
$d_t=-U_t-V_t$, and fix $R>0$ such that
$|t-d_t(L)|<R$ for all $L\in K$ and all sufficiently small $|t|$.
Define
\[
 I(t)=\frac{15}{4}\int_K\int_{[-R,R]^2}
 \bigl[U_t(L)+\min(0,X,Y,X+Y+t-d_t(L))\bigr]\,dX\,dY\,dL.
\]
The restriction of $I$ to $t<0$ has an analytic continuation $I_-$
across zero. For sufficiently small $t>0$,
\begin{equation}\label{eq:r4-common-area}
 \begin{split}
 I(t)-I_-(t)
 &=\frac58\int_K(t-d_t)_+^3\,dL\\
 &=\frac{15}{8}\int_0^t(t-s)^2B_t(s)\,ds,
 \qquad B_t(s)=\operatorname{area}\{L\in K:d_t(L)\le s\}.
 \end{split}
\end{equation}
Here $x_+=\max(x,0)$.
\end{lemma}
\begin{proof}
For fixed $L$, the required integral over the square is determined by
\begin{equation}\label{eq:r4-square-envelope}
 H(\delta)-H_-(\delta)=\frac{(\delta_+)^3}{6},\qquad
 H(\delta)=\int_{[-R,R]^2}\min(0,X,Y,X+Y+\delta)\,dX\,dY,
\end{equation}
where $|\delta|<R$ and $H_-$ is the polynomial continuation from
$\delta<0$. Indeed, $H'$ is the area on which the last plane is minimal.
For $\delta>0$ this is the rectangle $X,Y\le-\delta$; for $\delta<0$
one removes its corner triangle of area $\delta^2/2$. Integration from
zero proves \eqref{eq:r4-square-envelope}.

Now substitute $\delta=t-d_t(L)$. The hypotheses on the subdivision
and the affine coefficients make
\[
 I_-(t)=\frac{15}{4}\int_K
       \bigl[4R^2U_t(L)+H_-(t-d_t(L))\bigr]\,dL
\]
analytic near zero. Since $d_t\ge0$, the positive-part term vanishes
for $t<0$, so this expression agrees with $I(t)$ there. For $t>0$,
\eqref{eq:r4-square-envelope} gives the first equality in
\eqref{eq:r4-common-area}; integration over the sublevel sets of $d_t$
gives the second.
\end{proof}

To apply the lemma to the minimum
\[
 \min\{U_t(L),x,y,x+y+t+V_t(L)\},
\]
set $X=x-U_t(L)$ and $Y=y-U_t(L)$. This gives the integrand in $I(t)$.
For each $L$ the change is a translation, so it preserves measure by
Fubini's theorem, including where the affine expressions for $U_t$
change. We choose $K$ to contain the support of $(t-d_t)_+$. The
preceding localization then identifies the change $\Delta N(t)$ in the moment with
$I(t)-I_-(t)$. The two calculations below verify these conditions and
find the corresponding area $B_t(s)$.

\paragraph{The correction on $C_2$}

Put $t=2a$ and fix the other nine $C_2$ parameters. To distinguish the
edge parameter from the dual-coordinate vector, write $r_0=\lambda_{34}$.
Set
\[
 x_0=q_{11}-r_1,\quad y_0=q_{22}-r_2,\quad
 L_1=q_{33}+r_3,\quad L_2=q_{44}+r_4,\quad
 L_3=L_1+L_2-2r_0,
\]
and $(\alpha_1,\alpha_2,\alpha_3)=(p+r_0,q+r_0,p+q)$.
The change of four coordinates has absolute Jacobian one. The ten
site costs are
\[
 0,\ x_0,\ y_0,\ x_0+y_0+t,\quad
 L_j,\ x_0+y_0+t+2\alpha_j-L_j\quad(1\le j\le3).
\]
With
\[
 U=\min(0,L_1,L_2,L_3),\quad
 V=\min(0,2\alpha_1-L_1,2\alpha_2-L_2,2\alpha_3-L_3),\quad
 d=-U-V\ge0,
\]
their minimum becomes
$U+\min(0,X,Y,X+Y+t-d)$ after $X=x_0-U$, $Y=y_0-U$.
Here $U,V,d$ are independent of $t$.

Take $K_\eta=\{d\le\eta\}$, $R>2\eta$, and $|t|<\eta/4$,
all sufficiently small for the preceding localization. The functions
$U,V$ and their common subdivision on $K_\eta$ are independent of $t$.
Moreover, $|t-d|<R$ there, and $(t-d)_+^3$ is supported inside
$K_\eta$ for small positive $t$. Thus Lemma~\ref{lem:d4-area-integral}
applies.
The boundary $d=\eta$ has effective parameter $t-\eta<0$; the other
boundaries of the neighborhood stay away from the change of diagonal,
as required above.

The set $d\le s$ is given by
\[
\begin{gathered}
 -s\le L_1\le2(p+r_0)+s,\quad
 -s\le L_2\le2(q+r_0)+s,\\
 2r_0-s\le L_1+L_2\le2(p+q+r_0)+s,\\
 L_1-L_2\le2(p+r_0)+s,\quad L_2-L_1\le2(q+r_0)+s.
\end{gathered}
\]
These follow by expanding the sixteen pair inequalities in $-U-V\le s$;
the remaining eight are redundant. From the rectangle in the first line,
the next two lines remove two corner triangles of leg $2r_0+s$ and two
of leg $s$. Their interiors are disjoint, giving
\[
 B(s)=4(p+r_0+s)(q+r_0+s)-(2r_0+s)^2-s^2
     =4(pq+pr_0+qr_0)+4(p+q+r_0)s+2s^2.
\]
Substituting this area into \eqref{eq:r4-common-area} gives
\begin{align*}
 \Delta N
 &=\frac{15}{8}\int_0^{2a}(2a-s)^2B(s)\,ds\\
 &=2a^3\bigl[2a^2+5a(p+q+r_0)+10(pq+pr_0+qr_0)\bigr].
\end{align*}

\paragraph{The correction on $C_3$}
Put $t=2q_{34}$,
$u=q_{11}+q_{13}+q_{14}+t/2$,
$v=q_{22}+q_{23}+q_{24}+t/2$, and use coordinates
\[
 x=q_{33}-r_3,\quad y=q_{44}-r_4,\quad
 L_1=q_{11}+r_1,\quad L_2=q_{22}+r_2.
\]
The eight local costs have minimum
\[
 h_t=\min(0,x,y,x+y+t,L_1,L_2,
                      2u-L_1+x+y,2v-L_2+x+y).
\]
There is a path through actual quadratic forms with only the displayed
fourth plane depending on $t$: fix the diagonal entries and
$q_{12},q_{14},q_{24}$, and vary
$q_{34}=t/2$, $q_{13}=\alpha-t/2$, $q_{23}=\beta-t/2$.
Then $u,v$ and all coordinate offsets stay fixed. Write
$a_0,u_0,v_0,p_0,q_0$ for the values of $a,u,v,p,q$ at $b=0$, and
$j_0,k_0,l_0,m_0$ for those of $y_{13},y_{14},y_{23},y_{24}$,
respectively. In $C_3$ parameters the path is, with $b=t/2$,
\[
 (a_0-b,u_0,v_0,p_0-b,q_0+b,j_0+b,k_0,l_0+b,m_0,b).
\]
For small positive $b$ this lies in $C_3$; for small negative $b$ its
$\theta$ substitution lies in $C_2$.

For Lemma~\ref{lem:d4-area-integral}, use
\[
 U=\min(0,L_1,L_2),\qquad
 V_t=\min(0,A_1-L_1,A_2-L_2),\qquad
 A_1=2u-t,\quad A_2=2v-t.
\]
At the common boundary of $C_2$ and $C_3$, the dual face is the rectangle
$0\le L_1\le2u$, $0\le L_2\le2v$.
Take a fixed small neighborhood $K$ of this rectangle. The subdivisions
into regions where $U$ and $V_t$ are affine have vertices $(0,0)$ and $(A_1,A_2)$; their only
other crossings are $(A_1,0)$ and $(0,A_2)$. Because $A_1,A_2>0$
near the boundary, their common subdivision on a generic fixed
neighborhood has fixed combinatorial type and analytic vertices.
The affine coefficients depend linearly on $t$. For sufficiently small
positive $t$, the support of $(t-d_t)_+$ lies in $K$. Choosing this
neighborhood much thinner than the fiber square gives the uniform
bound $|t-d_t|<R$. These are the remaining conditions of
Lemma~\ref{lem:d4-area-integral}.

Expanding the nine pair inequalities in $-U-V_t\le s$ gives
\[
 -s\le L_1\le A_1+s,\quad -s\le L_2\le A_2+s,\quad
 L_1-L_2\le A_1+s,\quad L_2-L_1\le A_2+s.
\]
Thus its sublevel set is this rectangle with two opposite corner
triangles of leg $s$ removed, and
\[
 B_{2,t}(s)=(A_1+2s)(A_2+2s)-s^2
           =A_1A_2+2(A_1+A_2)s+3s^2.
\]
Substitution into the same integral, with $t=2b$, yields
\begin{align*}
 \Delta N
 &=\frac{15}{8}\int_0^t(t-s)^2B_{2,t}(s)\,ds\\
 &=2b^3\bigl[3b^2-5b(u+v)+10uv\bigr].
\end{align*}

\begin{proof}[Proof of Proposition~\ref{prop:d4-reconstructed-moments}]
The substitutions are formal polynomial substitutions and may have
negative coordinates. The determinant identities follow from
$A=E_{01}+E_{02}-E_{12}+E_{13}+E_{14}+E_{23}+E_{24}$ and
$E_b=A-E_{01}-E_{02}-E_{34}$.

The two calculations prove the moment identities near generic points
on the respective common boundaries, including open sets of the
tangential parameters. Fixed Delaunay
combinatorics makes the moment real analytic on each connected positive
parameter cone: dual vertices are linear in the matrix entries and
$\det Q>0$. For the boundary between $C_1$ and $C_2$, the negative-side expression
$5\Phi(z^{(2)})$ is an explicitly defined polynomial on the whole
parameter space. Its sum with the first correction therefore equals
the actual $N_2$ throughout $C_2$ by the identity theorem. This polynomial
then defines the formal expression $N_2(\theta)$, even where the
substituted coordinates are negative; the same argument proves the
second identity throughout $C_3$. Continuity covers positive-definite
boundary forms. Away from the neighborhood of the boundary between $C_2$ and $C_3$,
the second area integral denotes polynomial continuation and need not retain its geometric
positive-area interpretation.

Finally, expansion of the displayed graph polynomials gives nonnegative
coefficients in $N_i-20D(a+\frac14\sum_j\lambda_j)$: there are $367$
nonzero coefficients for $C_2$ and $373$ for $C_3$. Also
$[a^5]N_i=104$ and $\det A=4$, so $M(A)=104/(60\cdot4)=13/30$.
\end{proof}

\subsection{Polynomial inequalities}\label{subsec:polynomial-inequalities}
For $C_1$, order the ten parameters $x_{ij}$ by lexicographic order
on $0\le i<j\le4$ and put $L=\sum_{i<j}x_{ij}$.
For $C_2$ use the order
\[
 (a,x_{01},x_{02},x_{03},x_{04},
                       x_{13},x_{14},x_{23},x_{24},x_{34}),
\]
and for $C_3$ replace the final parameter $x_{34}$ by $b$.
In the latter two cases put
$t=a+\frac14\sum_{j=2}^{10}x_j$, with $x_1=a$.
Use the determinant polynomial $D$ and the moment polynomial $S$ on
$C_1$, or $N$ on $C_2,C_3$.

Let $H_1$ be the permutation group induced by $S_5$ on the ten edges of $K_5$.
For $C_2,C_3$ the permutation groups on the displayed parameter
positions $1,\ldots,10$ are
\begin{align*}
 H_2&=\langle
 (2\,6)(3\,8)(5\,10),
 (2\,3)(4\,5)(6\,9)(7\,8)\rangle,\\
 H_3&=\langle
 (2\,10)(4\,8)(5\,9)(6\,7),
 (2\,6\,10\,8\,3\,5)(4\,9\,7)\rangle.
\end{align*}
Their orders are $12$ and $36$, respectively.
These groups preserve the moment and determinant polynomials, and
preserve $L$ or $t$ as appropriate. For any polynomial $f$, write
\[
 \operatorname{Av}_{H} f=\frac1{|H|}\sum_{\sigma\in H}\sigma f.
\]

We also use the Schur cubic
\begin{equation}\label{eq:schur-cubic}
 \mathcal S_3(u,v,w)=u^3+v^3+w^3+3uvw
       -u^2v-u^2w-v^2u-v^2w-w^2u-w^2v\ge0
 \qquad(u,v,w\ge0).
\end{equation}
For completeness, by symmetry assume $u\ge v\ge w$. The polynomial
is $(u-v)^2(u+v-w)+w(u-w)(v-w)$, proving the inequality.

\begin{proposition}\label{prop:d4-exact-sos-certificates}
All six polynomials in the following table are nonnegative on the
closed nonnegative parameter orthant. The indicated appendix tables
specify integer polynomials $W_\nu$ with positive coefficients and
polynomials $B_\nu$ that are either squares $p_\nu^2$ or Schur cubics.
\[
\begin{array}{c|c|c|r|r|c}
\text{cone}&P&\beta&K&\text{rows}&\text{table}\\\hline
C_1&S-\tfrac43DL+\beta L^5&\tfrac7{4000}&1000&7&\ref{app:certificate-1-1}\\
C_1&S-\tfrac{3}{2}DL+\beta L^5&\tfrac{33}{10000}&1000&6&\ref{app:certificate-1-2}\\
C_2&N-\tfrac{82}{3}Dt+\beta t^5&\tfrac{41}{4}&100000&20&\ref{app:certificate-2-1}\\
C_2&N-31Dt+\beta t^5&20&10000&17&\ref{app:certificate-2-2}\\
C_3&N-\tfrac{82}{3}Dt+\beta t^5&\tfrac{41}{4}&100000&12&\ref{app:certificate-3-1}\\
C_3&N-31Dt+\beta t^5&20&100000&10&\ref{app:certificate-3-2}\\
\end{array}
\]
For each row, the explicitly defined remainder
\begin{equation}\label{eq:d4-sos-identity}
 R=P-\frac1K\sum_\nu\operatorname{Av}_H(W_\nu B_\nu)
\end{equation}
has nonnegative coefficients.
\end{proposition}

\begin{proof}
Every expression in \eqref{eq:d4-sos-identity} is a homogeneous
polynomial of degree five specified in the article: $D,S,N$ are the
graph formulas above, and all $W_\nu,B_\nu,K$ are displayed in
Appendix~\ref{app:explicit-certificates}. Direct expansion gives the
following coefficient comparisons, in the same order as the proposition:
\[
\begin{array}{c|r|r|c}
(\text{cone},k)&\#\{\alpha:[x^\alpha]R=0\}
 &\#\{\alpha:[x^\alpha]R>0\}&\min\{[x^\alpha]R:[x^\alpha]R>0\}\\\hline
(C_1,\tfrac43)&0&2002&\tfrac1{10000}\\
(C_1,\tfrac{3}{2})&50&1952&\tfrac1{60000}\\
(C_2,\tfrac{82}{3})&0&2002&\tfrac1{192000}\\
(C_2,31)&37&1965&\tfrac{73}{30000}\\
(C_3,\tfrac{82}{3})&0&2002&\tfrac{859}{7200000}\\
(C_3,31)&37&1965&\tfrac{41}{50000}\\
\end{array}
\]
Here $\alpha$ runs over the $\binom{14}{5}=2002$ nonnegative
integer exponent vectors of degree five. These comparisons can also
be made orbit by orbit: the groups $H_1,H_2,H_3$ have respectively
$35,225,91$ orbits on these monomials. For an orbit $\mathcal O$ and
any polynomial $f$,
\[
 [x^\alpha]\operatorname{Av}_H(f)
   =\frac1{|\mathcal O|}\sum_{\gamma\in\mathcal O}[x^\gamma]f
 \qquad(\alpha\in\mathcal O).
\]
In each case the zero and positive counts
sum to $2002$, so no coefficient is omitted. Since every $W_\nu$ has
positive coefficients and each $B_\nu$ is nonnegative by a square or
\eqref{eq:schur-cubic}, \eqref{eq:d4-sos-identity} proves $P\ge0$.
\end{proof}

The identities hold also on parameter faces. Continuity of the moment
therefore extends these inequalities to all positive-definite forms
in the closed cones.

\subsection{Completion of the global comparison}

The moment formulas and Proposition~\ref{prop:d4-exact-sos-certificates}
give, on $C_2$ and $C_3$,
\begin{equation}\label{eq:d4-assumed-bounds}
 N_i\ge\max\left\{20Dt,\;\frac{82}{3}Dt-\frac{41}{4}t^5,\;
                 31Dt-20t^5\right\},
 \qquad t=a+\frac14\sum_j\lambda_j.
\end{equation}
On $C_1$ they give
\begin{equation}\label{eq:d4-assumed-first}
 S\ge\max\left\{DL,\;\frac43DL-\frac7{4000}L^5,\;
                 \frac32DL-\frac{33}{10000}L^5\right\},
 \qquad L=\sum_{i<j}\lambda_{ij}.
\end{equation}

\begin{proof}[Proof of Theorem~\ref{thm:d4-global}]
For every rank-one generator $R_j$ of $C_2$ or $C_3$,
$\operatorname{tr}(A^{-1}R_j)=1$, while $\det A=4$.
The arithmetic--geometric mean inequality applied to
$A^{-1/2}QA^{-1/2}$ gives $D\le4t^4$, with equality precisely when
$Q=tA$. Set $d=D/(4t^4)\in(0,1]$. The required quantization bound is
equivalent to $N_i/t^5\ge104d^{5/4}$. By
\eqref{eq:d4-assumed-bounds}, the left side is at least
\[
 \max\{80d,\;\tfrac{328}{3}d-\tfrac{41}{4},\;124d-20\}.
\]
The three lines dominate $104d^{5/4}$ on
$(0,7/20]$, $[7/20,133/200]$, and $[133/200,1]$, respectively.
For the first interval this follows from $7/20<(10/13)^4$.
For each remaining interval, the difference between the line and
$104d^{5/4}$ is strictly concave, so it suffices to compare its endpoints.
At their endpoints the corresponding affine values are, respectively,
\[
 \frac{1681}{60},\quad \frac{18737}{300},\quad
 \frac{3123}{50},\quad104.
\]
Each value $v$ is positive and satisfies $v^4\ge104^4d^5$,
with equality only at the last endpoint $d=1$.
Thus the lower bound is strict for
$d<1$, and equality forces $Q=tA$.

On $C_1$, the matrix-tree theorem identifies $D$ with the weighted
spanning-tree polynomial of $K_5$. If $\rho_1,\ldots,\rho_4$ are the
nonzero eigenvalues of its weighted Laplacian, then
$\prod\rho_i=5D$ and $\sum\rho_i=2L$. Hence $D\le L^4/80$.
Scale $L$ to one and put $c_*=13/(120\sqrt2)$.
For $0<D\le1/200$, the first bound in
\eqref{eq:d4-assumed-first} gives
\[
 G(Q)\ge\frac1{48D^{1/4}}>c_*,
\]
since $(48c_*)^4=114244/625<200$.
For $D\in[1/200,9/1000]$ use the second bound in
\eqref{eq:d4-assumed-first}; for $D\in[9/1000,1/80]$ use the third.
On each interval, the difference between the corresponding affine
function and $48c_*D^{5/4}$ is strictly concave.
Its endpoint values are positive: if $\ell(D)$ denotes that affine
function, the required rational comparisons are
\[
 \ell(D)^4>\frac{114244}{625}D^5.
\]
They hold at both endpoints of each interval. Hence $G(Q)>c_*$
throughout $C_1$.
Since $M(A)=13/30$, the multiples of $A$ attain equality. To identify
the corresponding lattice, the matrix
\[
 B=\begin{pmatrix}1&1&-1&-1\\1&0&0&0\\0&1&0&0\\0&0&1&-1\end{pmatrix}
\]
satisfies $B^TB=A$ and $|\det B|=2$. Each column has even coordinate
sum, so $B\mathbb Z^4=D_4$. The extension to positive-definite boundary
forms and the classification by $C_1,C_2,C_3$ complete the proof.
\end{proof}

\section{Circular cones and lower bounds for lattice quantizers}
\label{sec:cone-bound}

For a bounded measurable set
$C\subset\mathbb R^d$ of positive volume, put
\[
 g(C)=\frac{\int_C\|x\|^2\,dx}
 {d\operatorname{vol}(C)^{1+2/d}}.
\]
The moment is taken about the specified origin; in particular,
$g(V(L))=G(L)$.

Put $A_d=d\kappa_d$. For $0<\omega<A_d/2$, let $g_d(\omega)$
denote $g(C)$ when $C$ is a right circular cone of solid angle
$\omega$ with vertex at the origin. Similarity invariance makes this
independent of the cone's size. A cone of height one and base radius
$r>0$ has
\begin{equation}\label{eq:intro-cone}
 \begin{split}
 \omega(r)&=(d-1)\kappa_{d-1}
       \int_0^r\frac{t^{d-2}}{(1+t^2)^{d/2}}\,dt,\\
 g_d(\omega(r))&=
 \frac{1+\dfrac{d-1}{d+1}r^2}
 {(d+2)\bigl(\kappa_{d-1}r^{d-1}/d\bigr)^{2/d}}.
 \end{split}
\end{equation}

We will use the convex function
\[
 F(v,j)=d^{d/2}\frac{v^{1+d/2}}{j^{d/2}},\qquad v,j>0.
\]
It is homogeneous of degree one. Moreover, $d^{-d/2}F(v,j)$ is the
perspective of the convex function $t\mapsto t^{1+d/2}$, so $F$ is convex.
Homogeneity and
Jensen's inequality give
\begin{equation}\label{eq:moment-subadditivity}
 F\left(\sum_i v_i,\sum_i j_i\right)\le\sum_i F(v_i,j_i).
\end{equation}
For a body $C$ of volume $v$ and second moment $j$, one has
$F(v,j)=g(C)^{-d/2}$.

\subsection{Concavity and the pyramid comparison}

\begin{lemma}\label{lem:cone-concavity}
For every $d\ge2$, the function
\[
 \phi_d(\omega)=g_d(\omega)^{-d/2},\qquad 0<\omega<A_d/2,
\]
is strictly concave.  It extends continuously by zero at both endpoints.
\end{lemma}
\begin{proof}
Scale the height to one, write $r$ for the radius of the base, and put
$\alpha=(d-1)/(d+1)$.  Direct integration gives the volume $v(r)$,
second moment $j(r)$, and solid angle $\omega(r)$ as
\begin{align}
 v(r)&=\frac{\kappa_{d-1}}d r^{d-1},&
 j(r)&=\frac{d}{d+2}v(r)(1+\alpha r^2),
 \label{eq:cone-volume-moment}\\
 \omega(r)&=(d-1)\kappa_{d-1}
       \int_0^{\arctan r}\sin^{d-2}\theta\,d\theta.
 \label{eq:cone-angle}
\end{align}
Consequently,
\[
 g_d(\omega(r))=\frac{1+\alpha r^2}{(d+2)v(r)^{2/d}},
 \qquad
 \phi_d(\omega(r))=K r^{d-1}(1+\alpha r^2)^{-d/2},
\]
where $K=(d+2)^{d/2}\kappa_{d-1}/d$.  If
$D=(d-1)\kappa_{d-1}$, then
\[
 \omega'(r)=D r^{d-2}(1+r^2)^{-d/2}>0,
 \qquad
 \phi_d'(\omega(r))=
 \frac KD\frac{[(d-1)-\alpha r^2](1+r^2)^{d/2}}
 {(1+\alpha r^2)^{d/2+1}}.
\]
Setting $s=r^2$, differentiation yields
\begin{equation}\label{eq:cone-derivative}
 \frac{d}{ds}
 \left(\frac{[(d-1)-\alpha s](1+s)^{d/2}}
 {(1+\alpha s)^{d/2+1}}\right)
 =-\frac{d(d-1)(d+2)}{(d+1)^2}
   \frac{s(1+s)^{d/2-1}}{(1+\alpha s)^{d/2+2}}<0.
\end{equation}
Both $s$ and $\omega(r)$ increase with $r$, proving strict concavity.
As $r\to0$ and $r\to\infty$, the solid angle tends to $0$ and $A_d/2$,
respectively.  The function $\phi_d(\omega(r))$ is asymptotic to a positive
multiple of $r^{d-1}$ at zero and of $r^{-1}$ at infinity, so both
endpoint limits vanish.
\end{proof}

The moment normalization in \cite{Ajjanagadde2020} is $d\,g(C)$.
Thus the function in its Conjecture~4 is
$d^{-d/2}\phi_d$, and Lemma~\ref{lem:cone-concavity} proves that
conjecture, with strict concavity.

\begin{lemma}\label{lem:cone-comparison}
Let $C\subset\mathbb R^d$, $d\ge2$, be the pyramid with vertex $0$ over
a compact convex base of dimension $d-1$ in an affine hyperplane not
containing $0$.  Suppose that the perpendicular foot from $0$ lies in
the relative interior of the base.  If $\omega$ is its solid angle at
$0$, then
\[
 g(C)\ge g_d(\omega).
\]
\end{lemma}
\begin{proof}
After scaling and rotating, write the base as
$\{(y,1):y\in B\}$, where $B\subset\mathbb R^{d-1}$ is a compact
convex body containing $0$ in its interior.  Let $\rho(u)>0$ be its
radial function on $S^{d-2}$, equipped with uniform probability
measure.  For $d=2$ this is the measure assigning weight $1/2$ to each
point of $S^0$.  Polar integration, with $v,j,\omega$ as in
\eqref{eq:cone-volume-moment}--\eqref{eq:cone-angle}, gives
\[
 \operatorname{vol}(C)=\mathbb E v(\rho),\qquad
 \int_C\|x\|^2\,dx=\mathbb E j(\rho),\qquad
 \omega=\mathbb E\omega(\rho).
\]
Indeed, each identity is obtained by integrating first along the rays
$\{tu:0\le t\le\rho(u)\}$ of the base.

Convexity of $F$, Jensen's inequality, and
Lemma~\ref{lem:cone-concavity} therefore imply
\begin{align*}
 g(C)^{-d/2}
 &=F(\mathbb E v(\rho),\mathbb E j(\rho))\\
 &\le\mathbb E F(v(\rho),j(\rho))
   =\mathbb E\phi_d(\omega(\rho))
   \le\phi_d(\mathbb E\omega(\rho))
   =\phi_d(\omega).
\end{align*}
Taking the power $-2/d$ proves the assertion.
\end{proof}

Lemma~\ref{lem:cone-comparison} proves the case of
Ajjanagadde's Conjecture~3 \cite[p.~139]{Ajjanagadde2020} in which the
perpendicular foot lies in the relative interior of the base. This
hypothesis holds for every pyramid over a lattice Voronoi facet, as
shown below.

\subsection{The lattice bound}

\begin{theorem}\label{thm:facet-bound}
Let $L\subset\mathbb R^d$ be a full-rank lattice, $d\ge2$, whose
Voronoi cell has $m$ facets.  Then
\begin{equation}\label{eq:actual-facet-lower}
 G(L)\ge m^{-2/d}g_d(A_d/m).
\end{equation}
In particular, with $f_d=2(2^d-1)$, every such lattice satisfies
\begin{equation}\label{eq:facet-lower}
 G(L)\ge f_d^{-2/d}g_d(A_d/f_d).
\end{equation}
\end{theorem}
\begin{proof}
Let $V$ be the Voronoi cell of $0$.
A facet separating $V$ and $V+\lambda$ is their intersection.  Since
$V=-V$, the map $x\mapsto\lambda-x$ preserves the facet, whose center
$\lambda/2$ is therefore in its relative interior.  This center is
also the perpendicular foot from $0$ to the supporting hyperplane.

Decompose $V$ into the pyramids over its $m$ facets, and write
$v_i,j_i,\omega_i$ for their volumes, second moments, and solid angles.
Their interiors are disjoint, and
\[
 \sum_i v_i=\operatorname{vol}V,\qquad \sum_i\omega_i=A_d.
\]
By \eqref{eq:moment-subadditivity},
Lemma~\ref{lem:cone-comparison}, and concavity,
\begin{align*}
 G(L)^{-d/2}
 &=F\left(\sum_i v_i,\sum_i j_i\right)\\
 &\le\sum_i F(v_i,j_i)
 \le\sum_i\phi_d(\omega_i)
 \le m\phi_d(A_d/m).
\end{align*}
Taking the power $-2/d$ proves \eqref{eq:actual-facet-lower}.

We recall the bound $m\le2(2^d-1)$.  If $\lambda$ defines a facet,
then $\pm\lambda$ are the only shortest vectors in its coset modulo
$2L$.  To see this, suppose that $w\equiv\lambda\pmod{2L}$,
$w\ne\pm\lambda$, and $\|w\|\le\|\lambda\|$.  The nonzero lattice
vectors
\[
 u=(\lambda+w)/2,\qquad z=(\lambda-w)/2
\]
satisfy $u+z=\lambda$.  Adding their Voronoi inequalities shows that
every $x\in V$ satisfies
\[
 2\langle x,\lambda\rangle
 \le\|u\|^2+\|z\|^2
 =\frac{\|\lambda\|^2+\|w\|^2}{2}.
\]
If $\|w\|<\|\lambda\|$, the supporting hyperplane for $\lambda$
does not meet $V$.  If equality holds, every point of its face must
satisfy both Voronoi inequalities for $u,z$ with equality.  In this
case $u$ and $z$ are nonzero and orthogonal, so the face has codimension
at least two.  Both alternatives contradict the facet assumption.
The zero coset has no facet vector, and each of the $2^d-1$ nonzero
cosets supplies at most one opposite pair.

Finally, concavity and $\phi_d(0)=0$ imply that $\phi_d(x)/x$ is
nonincreasing.  Since $m\le f_d$, we have
\[
 m\phi_d(A_d/m)\le f_d\phi_d(A_d/f_d),
\]
which proves \eqref{eq:facet-lower}.
\end{proof}

To compare \eqref{eq:facet-lower} with the original conjecture, let
$r_d$ be the unique positive solution of $\omega(r_d)=A_d/f_d$ and
put $h_d=(1+r_d^2)^{-1/2}$.  The cone with diagonal length one and
height $h_d$ in \cite[Conjecture~2]{Ajjanagadde2020} is similar to our
height-one cone of radius $r_d$.  Its moment $\xi(d,h_d)$ is
$d\,g_d(A_d/f_d)$, so the conjectured quantity
$\xi(d,h_d)/(d f_d^{2/d})$ is exactly the right-hand side of
\eqref{eq:facet-lower}.  Theorem~\ref{thm:facet-bound} therefore proves
the lattice assertion of that conjecture.

Let $R$ be the radial function of a right circular cone on its set of directions,
whose spherical measure is $\omega$. Its volume and second moment are
$\frac1d\int R^d$ and $\frac1{d+2}\int R^{d+2}$. Since $R$ is
nonconstant, strict Jensen's inequality gives
\[
 g_d(\omega)>\frac{1}{d+2}\left(\frac d\omega\right)^{2/d}.
\]
Taking $\omega=A_d/m$
shows that the lower bound in \eqref{eq:actual-facet-lower}, and in
particular that in \eqref{eq:facet-lower}, is strictly larger than
\eqref{eq:sphere-bound}.

For $d=2$, one has $f_2=6$, $A_2/f_2=\pi/3$, and
$r_2=1/\sqrt3$. Formula~\eqref{eq:cone-volume-moment} gives
\[
 f_2^{-1}g_2(\pi/3)=\frac{5}{36\sqrt3},
\]
the normalized second moment of the regular hexagon. Its six
triangles attain equality throughout the proof.

The Voronoi cell of $D_4$, by contrast, has $24$ regular octahedral
facets \cite{CS82}. The radial function of an octahedron about its
center is nonconstant. Strict concavity in the proof of
Lemma~\ref{lem:cone-comparison} therefore makes the pyramid comparison
strict, and hence
\[
 30^{-1/2}g_4(A_4/30)
 \le 24^{-1/2}g_4(A_4/24)<G(D_4).
\]

\appendix
\section{Explicit polynomial certificates}\label{app:explicit-certificates}
The letters $(a,b,c,d,e,f,g,h,i,j)$ denote, in order, the ten
parameters of Section~\ref{subsec:polynomial-inequalities}. Thus
$L=a+\cdots+j$ on $C_1$, while $t=a+(b+\cdots+j)/4$ on $C_2,C_3$.
The polynomials below give the identity
\[
 P=R+\frac1K\sum_\nu\operatorname{Av}_H(W_\nu B_\nu)
\]
of Proposition~\ref{prop:d4-exact-sos-certificates}. In an unmarked
row $B_\nu=p_\nu^2$; in a row marked $*$ the middle column gives
$(u,v,w)$ and $B_\nu=\mathcal S_3(u,v,w)$.
Every nonzero $W_\nu$ has positive coefficients. The remainder $R$
is specified by this identity, and its coefficient comparison is given
in the proposition. Table~\ref{app:certificate-1-1} uses one list of
linear forms for both $C_1$ inequalities; a zero weight contributes
nothing.

\newcommand{\certrow}[3]{#1&$#2$&$#3$\\[4pt]}
\newcommand{\certpair}[4]{#1&$#2$&$#3$&$#4$\\[4pt]}
\begingroup
\footnotesize
\setlength{\tabcolsep}{3pt}
\renewcommand{\arraystretch}{1.08}
Here $K=1000$ and $H=H_1$. The weights $W^{(1)}$ and $W^{(2)}$
correspond to $k=4/3$ and $k=3/2$, respectively.
\begin{longtable}{@{}r p{0.29\linewidth} p{0.32\linewidth} p{0.32\linewidth}@{}}
\caption{The two $C_1$ certificates, with $K=1000$.}\label{app:certificate-1-1}\label{app:certificate-1-2}\\
\toprule
$\nu$ & $p_\nu$ & $W^{(1)}_\nu$ & $W^{(2)}_\nu$\\\midrule
\endfirsthead
\toprule
$\nu$ & $p_\nu$ & $W^{(1)}_\nu$ & $W^{(2)}_\nu$\\\midrule
\endhead
\bottomrule
\endfoot
\certpair{1}{\begin{aligned}[t]&a+b+c-(d+h)\end{aligned}}{\begin{aligned}[t]&16bd^{2}+159efh\end{aligned}}{\begin{aligned}[t]&43d^{2}h+3737deh+4dh^{2}+368efh\end{aligned}}
\certpair{2}{\begin{aligned}[t]&a+b+c-(h+i+j)\end{aligned}}{\begin{aligned}[t]&30b^{2}h+93bfh+64egh+364gij\\
&+16h^{3}\end{aligned}}{\begin{aligned}[t]&214abh+226bh^{2}+102ehi\\
&+122gh^{2}+2192ghi+11h^{3}+127h^{2}i\\
&+175hi^{2}\end{aligned}}
\certpair{3}{\begin{aligned}[t]&a+b+c-(g+i+j)\end{aligned}}{\begin{aligned}[t]&1716aef+58aeh+953e^{2}f\\
&+120efh\end{aligned}}{\begin{aligned}[t]&327a^{2}b+411a^{2}h+200abc\\
&+626abe+405abf+201abg+43ae^{2}\\
&+850aef+2065aeh+266agh\\
&+687ah^{2}+1284e^{2}f+175efh\end{aligned}}
\certpair{4}{\begin{aligned}[t]&a+b+c+d-(i+j)\end{aligned}}{\begin{aligned}[t]&223aij+699ehj+2580eij\\
&+121ej^{2}+400gij+263i^{2}j\end{aligned}}{0}
\certpair{5}{\begin{aligned}[t]&a+b+c+d-(g+h)\end{aligned}}{\begin{aligned}[t]&96a^{2}g+947egh+47g^{2}h\end{aligned}}{\begin{aligned}[t]&637egh\end{aligned}}
\certpair{6}{\begin{aligned}[t]&a+b+e-(f+g+h+i)\end{aligned}}{\begin{aligned}[t]&2abf+537abj+110ae^{2}+23aej\\
&+30fgj\end{aligned}}{\begin{aligned}[t]&699abj+656aej+22ej^{2}+41j^{3}\end{aligned}}
\certpair{7}{\begin{aligned}[t]&a+b+c+d-(h+i+j)\end{aligned}}{0}{\begin{aligned}[t]&21a^{3}+59a^{2}h+772ehi+1275ehj\end{aligned}}
\certpair{8}{\begin{aligned}[t]&L-2(d+h)-(e+f)\end{aligned}}{\begin{aligned}[t]&1979deh+121dh^{2}\end{aligned}}{0}
\end{longtable}
\begin{longtable}{@{}r p{0.37\linewidth} p{0.56\linewidth}@{}}
\caption[Certificate on C2 with k=82/3]{The $C_2$ certificate with $k=\tfrac{82}{3}$, $K=100000$, and $H=H_2$.}\label{app:certificate-2-1}\\
\toprule
$\nu$ & $p_\nu$ & $W_\nu$\\\midrule
\endfirsthead
\toprule
$\nu$ & $p_\nu$ & $W_\nu$\\\midrule
\endhead
\bottomrule
\endfoot
\certrow{1}{\begin{aligned}[t]&h+j-i\end{aligned}}{\begin{aligned}[t]&105003b^{2}i+29522bci+3768bd^{2}+26153bdi\\
&+1814bef+14569bei+3898bf^{2}+57180bgi+32485bi^{2}\\
&+15791efi+8565fj^{2}\end{aligned}}
\certrow{2}{\begin{aligned}[t]&b+j-(c+h+i)\end{aligned}}{\begin{aligned}[t]&51553abj+126524bdj+162218bfj+25700cdj\\
&+6366dej+6053dj^{2}+22862fgj\end{aligned}}
\certrow{3}{\begin{aligned}[t]&10(b+f)-9(c+h)-5i\end{aligned}}{\begin{aligned}[t]&407abf+330ad^{2}+645ae^{2}+120b^{2}f+3bcd+440bde\\
&+421bdf+1034bdg+2573bdj+502bef+243beg\\
&+276bfg+741cd^{2}+151d^{3}+1106d^{2}e+181d^{2}g+121de^{2}\\
&+49dej+565e^{2}j+69ei^{2}+22i^{3}\end{aligned}}
\certrow{4}{\begin{aligned}[t]&b+c+d-(g+i+j)\end{aligned}}{\begin{aligned}[t]&59413adf+6579bd^{2}+2921f^{3}+63978f^{2}h\end{aligned}}
\certrow{5}{\begin{aligned}[t]&4t-4a-2(f+g+j)-c\end{aligned}}{\begin{aligned}[t]&165234afg+334071afj+4815bfj+7456cde\\
&+1199cf^{2}+93229cfg+195424cfj+228dfg+23115dfj\\
&+10901fij\end{aligned}}
\certrow{6}{\begin{aligned}[t]&(c-j)^{2}\end{aligned}}{\begin{aligned}[t]&2638a+20229b+1995f+2895j\end{aligned}}
\certrow{7}{\begin{aligned}[t]&3a+c-j\end{aligned}}{\begin{aligned}[t]&46372a^{2}j+69567adj+146590afj+27687bcd\\
&+14531bhj+92085dfj+143024dgj+13298fhj\\
&+29056fij\end{aligned}}
\certrow{8}{\begin{aligned}[t]&2(c+e+h+j)-3i\end{aligned}}{\begin{aligned}[t]&1975bd^{2}+4101bdg+7340d^{2}i+190dgi+4066di^{2}\end{aligned}}
\certrow{9}{\begin{aligned}[t]&5i+4c-5f-4(d+j)\end{aligned}}{\begin{aligned}[t]&1379abc+32ai^{2}+845b^{2}c+7676b^{2}e+5712b^{2}g\\
&+675bce+2222bcg+4590be^{2}+4604beg+300bf^{2}\\
&+1106bfg+1222bgj+883bj^{2}+537cei+183dej\\
&+2661dfg+273ej^{2}+1456g^{2}j+382gj^{2}\end{aligned}}
\certrow{10}{\begin{aligned}[t]&2a+c+h+i-2g\end{aligned}}{\begin{aligned}[t]&123391a^{2}g+117038abg+16868acg+66137adg\\
&+104443aeg+22003agi+16327b^{2}g+29689bdg\\
&+9048bg^{2}+20852bgj+8672cd^{2}+1144cdi+1629ce^{2}\\
&+1985d^{3}+4829dgi+31248egj\end{aligned}}
\certrow{11}{\begin{aligned}[t]&c+h+i-(d+e+j)\end{aligned}}{\begin{aligned}[t]&21209ab^{2}+17628abc+12014abj+13152acj\\
&+117861ade+77608bde+124279bdj+1206bfg\\
&+18769d^{2}e+488dej\end{aligned}}
\certrow{12}{\begin{aligned}[t]&c+e+g+h+j-(f+i)\end{aligned}}{\begin{aligned}[t]&22415adi+44786ai^{2}+3539bde+189262bfi\\
&+251709dfi\end{aligned}}
\certrow{13}{\begin{aligned}[t]&c+e+i-(d+f+h+j)\end{aligned}}{\begin{aligned}[t]&71686ac^{2}+100179ace+69769aci+770b^{2}c+8270b^{2}g\\
&+130172bce+178859bci+5632bd^{2}+6627be^{2}\\
&+114978bei+37153bi^{2}+4428c^{3}+29402c^{2}e\\
&+13136chi\end{aligned}}
\certrow{14}{\begin{aligned}[t]&3d+2(f+h+j)+a-3i-2g\end{aligned}}{\begin{aligned}[t]&3412abi+3256afi+29594agi+143bci+12738bgi\\
&+9148cgi+1250egi+2168g^{2}i+5367gij+739i^{2}j\end{aligned}}
\certrow{15}{\begin{aligned}[t]&4t-2(g+i)-(d+e)\end{aligned}}{\begin{aligned}[t]&19076agi+3537d^{3}+95798dgi+5346fgi\end{aligned}}
\certrow{16}{\begin{aligned}[t]&2(a+g)+c+e+h+j-2d-i\end{aligned}}{\begin{aligned}[t]&149666adi+37772bdi+1759d^{2}i+1732gi^{2}\end{aligned}}
\certrow{17}{\begin{aligned}[t]&4t-3e-(a+f+h)\end{aligned}}{\begin{aligned}[t]&3461ef^{2}+46017efh\end{aligned}}
\certrow{18}{\begin{aligned}[t]&8t-8a-5(f+i)\end{aligned}}{\begin{aligned}[t]&2165adi+989afh+37075afi+1129fgi\end{aligned}}
\certrow{19$^*$}{(b,c,d)}{\begin{aligned}[t]&21852bd+631de+2258df+9515eg+8770fg+46347fi\end{aligned}}
\certrow{20$^*$}{(b,c,j)}{\begin{aligned}[t]&233579a^{2}+146665ad+191173af+57313aj+820b^{2}\end{aligned}}
\end{longtable}
\begin{longtable}{@{}r p{0.37\linewidth} p{0.56\linewidth}@{}}
\caption[Certificate on C2 with k=31]{The $C_2$ certificate with $k=31$, $K=10000$, and $H=H_2$.}\label{app:certificate-2-2}\\
\toprule
$\nu$ & $p_\nu$ & $W_\nu$\\\midrule
\endfirsthead
\toprule
$\nu$ & $p_\nu$ & $W_\nu$\\\midrule
\endhead
\bottomrule
\endfoot
\certrow{1}{\begin{aligned}[t]&3a-2j\end{aligned}}{\begin{aligned}[t]&1433ab^{2}+3267abj+93b^{2}c+8321bcj+1828bfj\\
&+265fhj\end{aligned}}
\certrow{2}{\begin{aligned}[t]&2c-3j\end{aligned}}{\begin{aligned}[t]&8555a^{2}j+6083abj+232ac^{2}+4145acd+2162adj\\
&+12042afj+1285aj^{2}+1826bc^{2}+1330bde+1872bdj\\
&+561bfj+318d^{2}j+358dgj+256dj^{2}+43fj^{2}+234h^{2}j\\
&+1441hij\end{aligned}}
\certrow{3}{\begin{aligned}[t]&g-i\end{aligned}}{\begin{aligned}[t]&15975a^{3}+6024a^{2}g+8759abg+4684ad^{2}+114568adg\\
&+1128ag^{2}+12186agi+12832bd^{2}+18070bdi\\
&+20334d^{2}e+1559egi+872g^{3}\end{aligned}}
\certrow{4}{\begin{aligned}[t]&e+i-(f+j)\end{aligned}}{\begin{aligned}[t]&9055a^{2}f+1958b^{2}c+33214bc^{2}+6711c^{3}+2067c^{2}g\\
&+3572c^{2}j+3314dg^{2}+6367dgi\end{aligned}}
\certrow{5}{\begin{aligned}[t]&c+i-(d+j)\end{aligned}}{\begin{aligned}[t]&11213a^{3}+6196a^{2}d+20429a^{2}e+17942ace+11608ade\\
&+28412adj+5841b^{2}g+3821bde+25811bdj+181c^{2}e\end{aligned}}
\certrow{6}{\begin{aligned}[t]&c+d-(i+j)\end{aligned}}{\begin{aligned}[t]&11213a^{3}+31215a^{2}h+28355ab^{2}+2227abg+21035abh\\
&+17340acd+41447ach+23561adh+14948ah^{2}\\
&+23821b^{2}h+7938b^{2}i+6634bcd+10461bdh+11339bhj\\
&+10642bij+3744ch^{2}+5801chi+18714dh^{2}+1541dhj\\
&+1394h^{3}\end{aligned}}
\certrow{7}{\begin{aligned}[t]&f+g-(h+i)\end{aligned}}{\begin{aligned}[t]&24137a^{2}j+4799abc+6895abh+27623abj+10089afg\\
&+62418afj+30001aj^{2}+249b^{3}+2527b^{2}h+3142bfj\\
&+5451bhi+20299bhj+16115bj^{2}+6796dfj+171dj^{2}\\
&+2080f^{2}g+8964f^{2}j+1692fgj+8786fij+7603fj^{2}\\
&+7j^{3}\end{aligned}}
\certrow{8}{\begin{aligned}[t]&2(d+f)+h+j-4i\end{aligned}}{\begin{aligned}[t]&2076a^{2}b+3128a^{2}i+5803abf+7205abi+3813aei\\
&+1221aij+159b^{2}d+1818b^{2}h+373bcd+2665bci\\
&+631bd^{2}+2529bde+1269bdi+2507bei+273bf^{2}\\
&+652bfi+3064bgi+877bhi+127bij+1107cd^{2}\\
&+3304cdi+549cgi+109d^{3}+298d^{2}f+1826d^{2}g+238d^{2}h\\
&+1875dei+1526dg^{2}+1876dgi+628ehi+557fi^{2}\\
&+659g^{2}i+1228gij+668h^{2}i+356hi^{2}+123i^{2}j+331ij^{2}\end{aligned}}
\certrow{9}{\begin{aligned}[t]&2(f+i)-2(g+h)-j\end{aligned}}{\begin{aligned}[t]&14992a^{2}f+2140abf+7525afi+923bdf+2104bdi\\
&+3053bfi+2217dfi+552f^{2}i\end{aligned}}
\certrow{10}{\begin{aligned}[t]&10i+9c-10f-7(d+j)\end{aligned}}{\begin{aligned}[t]&154a^{2}i+365ab^{2}+20abe+378abg+136ae^{2}+114aej\\
&+423b^{2}e+502b^{2}g+120bci+175bde+272be^{2}+51bef\\
&+231beg+152bei+16bf^{2}+138bfg+221bfj+126bgj\\
&+100dfj+57dj^{2}+47e^{2}f+126e^{2}g+54e^{2}h+56ef^{2}\\
&+84efh+348egi\end{aligned}}
\certrow{11}{\begin{aligned}[t]&2(d+j)+b+g-2(c+i)\end{aligned}}{\begin{aligned}[t]&2921a^{2}e+984a^{2}f+6816acf+5703aci+356adf\\
&+12585aef+1112bfj+4465cef+850cf^{2}+3450cfi\\
&+4254def+355e^{2}f+4122ef^{2}\end{aligned}}
\certrow{12}{\begin{aligned}[t]&6b+5f+4g-6c-5h-4i\end{aligned}}{\begin{aligned}[t]&444a^{2}b+535a^{2}d+161ab^{2}+80abj+868ade+262bde\\
&+246bdg+457bdj+128d^{2}e+311de^{2}+74dej\end{aligned}}
\certrow{13}{\begin{aligned}[t]&2g+c+e+h+j-2d-i\end{aligned}}{\begin{aligned}[t]&2602a^{2}d+6173abd+10265adi+7229bdi+1761cdi\end{aligned}}
\certrow{14}{\begin{aligned}[t]&4d+3f+2(h+j)+a-4i-2g\end{aligned}}{\begin{aligned}[t]&1801agi+576bcd\end{aligned}}
\certrow{15}{\begin{aligned}[t]&4t-4a-2(f+g+j)-c\end{aligned}}{\begin{aligned}[t]&4953a^{2}c\end{aligned}}
\certrow{16}{\begin{aligned}[t]&3(c+d+e)+2b-3f-2(g+j)-i\end{aligned}}{\begin{aligned}[t]&3414a^{2}f+678acf+2311aci+1675af^{2}+5111afh\\
&+207afi+224ai^{2}+2572aij+1049bfh+185cfh\\
&+2080fgh+9fgi\end{aligned}}
\certrow{17}{\begin{aligned}[t]&4t+b+c+f-3(g+i)-d\end{aligned}}{\begin{aligned}[t]&31dgi\end{aligned}}
\end{longtable}
\begin{longtable}{@{}r p{0.37\linewidth} p{0.56\linewidth}@{}}
\caption[Certificate on C3 with k=82/3]{The $C_3$ certificate with $k=\tfrac{82}{3}$, $K=100000$, and $H=H_3$.}\label{app:certificate-3-1}\\
\toprule
$\nu$ & $p_\nu$ & $W_\nu$\\\midrule
\endfirsthead
\toprule
$\nu$ & $p_\nu$ & $W_\nu$\\\midrule
\endhead
\bottomrule
\endfoot
\certrow{1}{\begin{aligned}[t]&2(e+i)+b-2(d+h)-j\end{aligned}}{\begin{aligned}[t]&9025aef+10444af^{2}+28552b^{2}g+16596bc^{2}\\
&+29113bcg+4493bci+4579bdh+19671bfj+17003bg^{2}\\
&+670bgh+10659c^{3}+25634c^{2}d+50724c^{2}g+7572cde\\
&+50317cdg+8453cdi+4913cg^{2}\end{aligned}}
\certrow{2}{\begin{aligned}[t]&d+f+h+j-(g+i)\end{aligned}}{\begin{aligned}[t]&221608agi+383059bgi+141120egi\end{aligned}}
\certrow{3}{\begin{aligned}[t]&3b+2h+j-3(c+i)-2(d+f)\end{aligned}}{\begin{aligned}[t]&32508abj+8096bdg+5228be^{2}+60756beh+25911bgh\\
&+30745bhj+2331dhi+1257h^{2}j\end{aligned}}
\certrow{4}{\begin{aligned}[t]&4t+g-4a-3(c+d+f+i)\end{aligned}}{\begin{aligned}[t]&6483abg+2293acg+37759b^{2}j+4685bdg+143714bgh\\
&+136167bgj\end{aligned}}
\certrow{5}{\begin{aligned}[t]&12t-10a-8(e+i)-3c\end{aligned}}{\begin{aligned}[t]&8957a^{2}e+15670ace+2269aef+39281aei\end{aligned}}
\certrow{6}{\begin{aligned}[t]&4t+6(a+g)+c+f-6(d+i)\end{aligned}}{\begin{aligned}[t]&5583a^{2}d+14250adi+482bdh+995cdf+452d^{2}i\\
&+1267dgi\end{aligned}}
\certrow{7}{\begin{aligned}[t]&12t-5(b+c+j)-2a\end{aligned}}{\begin{aligned}[t]&16836bcj\end{aligned}}
\certrow{8}{\begin{aligned}[t]&12t-11a-5(d+e+j)\end{aligned}}{\begin{aligned}[t]&6403ade+13049dej\end{aligned}}
\certrow{9}{\begin{aligned}[t]&20t-17a-13(h+i)\end{aligned}}{\begin{aligned}[t]&2150a^{2}h+673ah^{2}+8042ahi+34b^{2}h+80bh^{2}+858chi\\
&+878h^{2}i\end{aligned}}
\certrow{10$^*$}{(b,c,f)}{\begin{aligned}[t]&121684a^{2}+4624bc+15384ch+62700eh+48753ei\\
&+103934fi+16225h^{2}+9591hi\end{aligned}}
\certrow{11$^*$}{(b,c,j)}{\begin{aligned}[t]&202461a^{2}+136899ab+354697ad+6136bc\end{aligned}}
\certrow{12$^*$}{(b,d,f)}{\begin{aligned}[t]&8831b^{2}+19145bc+44233be+27608bh+76500cg\end{aligned}}
\end{longtable}
\begin{longtable}{@{}r p{0.37\linewidth} p{0.56\linewidth}@{}}
\caption[Certificate on C3 with k=31]{The $C_3$ certificate with $k=31$, $K=100000$, and $H=H_3$.}\label{app:certificate-3-2}\\
\toprule
$\nu$ & $p_\nu$ & $W_\nu$\\\midrule
\endfirsthead
\toprule
$\nu$ & $p_\nu$ & $W_\nu$\\\midrule
\endhead
\bottomrule
\endfoot
\certrow{1}{\begin{aligned}[t]&c-j\end{aligned}}{\begin{aligned}[t]&457158a^{3}+640641ab^{2}+365125b^{2}c+683333b^{2}d\\
&+313473b^{2}h+169301bcf+501569bch+142614cde\end{aligned}}
\certrow{2}{\begin{aligned}[t]&2a-j\end{aligned}}{\begin{aligned}[t]&1940ab^{2}+106768aj^{2}+77368bfj+75261fhj\end{aligned}}
\certrow{3}{\begin{aligned}[t]&b+i-(h+j)\end{aligned}}{\begin{aligned}[t]&147029ae^{2}+218144aij+21498b^{2}c+122288b^{2}j\\
&+97510bch+88472bci+663082bei+632030bgi\\
&+141649bij+23170cij+185330egi+4425gij\\
&+89635hij+124366i^{2}j\end{aligned}}
\certrow{4}{\begin{aligned}[t]&2d+f+h+j-3i-e\end{aligned}}{\begin{aligned}[t]&82621a^{2}b+1751a^{2}h+350478a^{2}i+10726ab^{2}\\
&+21772abh+42314abi+102719aci+12326afi\\
&+302203agi+90977aij+2666bci+49601egi\end{aligned}}
\certrow{5}{\begin{aligned}[t]&2(b+i)-(c+d+g+j)\end{aligned}}{\begin{aligned}[t]&163896a^{2}h+201338abi+11325ah^{2}+85bci+3051bh^{2}\\
&+78453bhi\end{aligned}}
\certrow{6}{\begin{aligned}[t]&b+g+j-(c+d+i)\end{aligned}}{\begin{aligned}[t]&82668a^{2}c+857862a^{2}e+519589abd+23379abg\\
&+706867abh+398759ace+90817ae^{2}+945256aeh\\
&+22884b^{2}d+139028b^{2}h+41470bcd+76577beg\\
&+259597beh+351713bej+245213bgh+198142bgj\\
&+309245bh^{2}+121695ceh+12932e^{3}+122877e^{2}h\end{aligned}}
\certrow{7}{\begin{aligned}[t]&3e+2(b+i)-3h-2(d+j)\end{aligned}}{\begin{aligned}[t]&6158abc+72386abe+82420ac^{2}+21648bc^{2}\\
&+26048bce+77329bcg+8292c^{3}+9964c^{2}e+21090c^{2}g\\
&+29182ceg+2593e^{2}h\end{aligned}}
\certrow{8}{\begin{aligned}[t]&4t-3(d+g+i)-a\end{aligned}}{\begin{aligned}[t]&1496adg+4100bcj+8214bdg+38215dgi\end{aligned}}
\certrow{9}{\begin{aligned}[t]&3(b-h)(d-g)+b(f-e)+f(d+f+i)\\
&+j(g+h+i+j)-c(c+d+h+i)\\
&-e(e+g+i)\end{aligned}}{\begin{aligned}[t]&46259a+12227c+35068i\end{aligned}}
\certrow{10}{\begin{aligned}[t]&2(b+j)^{2}+2bj\\
&+(b+j)(f+g-2d-2e-2h-2i)\\
&-3c(f+g)-3(f^{2}+fg+g^{2})\\
&+(4f+g)(d+i)+(f+4g)(e+h)\\
&+4(di+eh)-2(de+hi)\end{aligned}}{\begin{aligned}[t]&2266a+81d+500f\end{aligned}}
\end{longtable}
\endgroup

\section*{Acknowledgments}
Generative AI tools were used to discuss research ideas, construct and
verify the polynomial certificates, assist in checking the correctness
of the proofs, and proofread the manuscript.
The author is responsible for the correctness and content of the paper.

\begingroup
\renewcommand{\bibliofont}{\fontsize{8}{9}\selectfont}
\raggedright
\raggedbottom
\bibliographystyle{alpha}
\bibliography{references}
\endgroup

\end{document}